\documentclass[10pt,english]{amsart}

\usepackage[T1]{fontenc}
\usepackage[latin9]{inputenc}
\usepackage{amsmath}
\usepackage{amsthm}
\usepackage{txfonts}
\usepackage[letterpaper,margin=3cm]{geometry}
\usepackage{babel}

\theoremstyle{plain}

\newtheorem{theorem}{Theorem}[section]

\newtheorem{lemma}[theorem]{Lemma}
\newtheorem{proposition}[theorem]{Proposition}

\newcommand{\RealVect}[1]{{\mathbb R}^{#1}}

\newcommand{\cws}{\stackrel{*}{\to}} 
\newcommand{\Ud}[1]{\tilde{U}_d^{#1}} 
\newcommand{\Rk}[3]{\frac{1}{|{#1} - {#2}|^{#3}}} 
\newcommand{\Hdrest}[1]{\Hd_{#1}} 

\def\Rl{\mathbb R} 
\def\Rp{\RealVect{p}} 
\def\MA{{\mathcal M}(A)} 
\def\MAp{{\mathcal M}_1(A)}
\def\Es{{\mathcal E}_s} 
\def\Hd{{\mathcal H}^d} 
\def\Hdr{\Hdrest{A}} 
\def\Id{\tilde{I}_d} 
\def\Is{I_s} 

\DeclareMathOperator{\diam}{diam}
\DeclareMathOperator{\dist}{dist}
\DeclareMathOperator{\supp}{supp}

\begin{document}

\title{Riesz $s$-equilibrium measures on $d$-dimensional fractal sets as $s$ approaches $d$}

\author{Matthew T. Calef}
\address{M. T. Calef:
Department of Mathematics, 
Vanderbilt University, 
Nashville, TN 37240, 
USA }
\email{Matthew.T.Calef@Vanderbilt.Edu}

\begin{abstract}
Let $A$ be a compact  set in $\Rp$ of Hausdorff dimension $d$.  For $s\in(0,d)$, the Riesz $s$-equilibrium measure $\mu^{s,A}$ is the unique Borel probability measure   with support in $A$ that
minimizes  $$ \Is(\mu):=\iint\Rk{x}{y}{s}d\mu(y)d\mu(x)$$  over all such probability measures.   In this paper we show that if  $A$ is a strictly self-similar $d$-fractal, then $\mu^{s,A}$ converges in the weak-star topology  to normalized $d$-dimensional Hausdorff measure restricted to $A$  as $s$ approaches $d$ from below.  
\end{abstract}

\maketitle


\section{Introduction}

Let $A$ be a compact subset of $\Rp$ with positive $d$-dimensional Hausdorff measure. Let $\MA$ denote the (unsigned) Radon measures supported on $A$ and $\MAp\subset\MA$ the probability measures in $\MA$. Recall (cf.~\cite{Landkof1, Fuglede1, Mattila1}) that for $s\in(0,d)$ the \emph{Riesz $s$-energy} of a measure $\mu\in\MA$ is $$I_s(\mu) := \iint \Rk{x}{y}{s} d\mu(x)d\mu(y),$$ and that there is a unique measure $\mu^{s,A}\in \MAp$ called the {\em equilibrium measure} with the property that $I_s(\mu^{s,A}) < I_s(\nu)$ for all $\nu \in \MAp \backslash \{\mu^{s,A}\}$. The \emph{$s$-potential} of a measure $\mu$ at a point $x$ is $$U_s^\mu(x) := \int \Rk{x}{y}{s}d\mu(y),$$ and for any measure $\mu$ with finite $s$-energy $$I_s(\mu) = \int U_s^\mu d\mu.$$ For $s \ge d$, $I_s(\mu) =\infty$ for all non-trivial $\mu\in\MA$. We shall denote the $d$-dimensional Hausdorff measure as $\Hd$ and the restriction of a measure $\mu$ to a set $E$ as $\mu_E$ e.g. $\Hdr := \Hd(\cdot \cap A)$. The closed ball of radius $r$ centered at $x$ is denoted $B(x,r)$. 

The study of equilibrium measures arises naturally in electrostatics. Specifically one may consider $\mu^{s,A}$ as the positive charge distribution on $A$ that minimizes a generalized electrostatic energy mediated by the kernel $|x-y|^{-s}$ where, in the classical electrostatic or Newtonian setting $s=d-2$. In the case of the interval $A=[-1,1]$ ($d=1$) it is known $\mu^{s,[-1,1]}$ is absolutely continuous with respect to the one-dimensional Lebesgue measure and the Radon-Nikod\'ym derivative of $\mu^{s,[-1,1]}$ is $c_s(1-x^2)^{\frac{1-s}{2}}$ where $c_s$ is chosen to make the measure of unit mass. One can see $\mu^{s,[-1,1]}$ converges in the weak-star topology on $\MA$ to $\Hdr/\Hd(A)$ as $s\uparrow 1$. More generally it is shown in~\cite{CalefHardin1} that this convergence occurs for certain $d$-rectifiable sets. 

In this paper we prove the same result for any compact self-similar fractal $A\subset \Rp$ satisfying $$A = \bigcup_{i=1}^N \varphi_i(A),$$ where the union is disjoint and the maps $\varphi_1,\ldots,\varphi_N$ satisfy $|\varphi_i(x)| = L_i|x|$ for all $x\in \Rp$ and where $L_i\in(0,1)$.  We refer to such sets as \emph{strictly self-similar $d$-fractals}. In \cite{Moran1} Moran shows for strictly self-similar $d$-fractals the Hausdorff dimension is also the unique value of $d$ that satisfies the equation $$\sum_{i=1}^N L_i^d = 1,$$ and that $\Hd(A)\in(0,\infty)$. Moran shows this results for fractals satisfying the broader \emph{open set condition} (cf.~\cite{Falconer1}), however we use the strict separation in the proofs of the following results.

Given a Borel measure $\mu$, let $\Theta_d^r(\mu,x):=\mu(B(x,r))/r^d$ denote the average $d$-density of $\mu$ over a radius $r$ about $x$. The limit as $r\downarrow 0$, $$\Theta_d(\mu,x) := \lim_{r\downarrow 0} \Theta_d^r(\mu,x),$$ when it exists, is the classical point density of $\mu$ at $x$. It is consequence of a result of Preiss~\cite{Preiss1} (also cf.~\cite{Mattila1}) that if $A$ is a strictly self-similar $d$-fractal, then at $\Hdr$-a.a. $x\in A$ the point densities $\Theta_d(\Hdr,x)$ do not exist. However, Bedford and Fisher in ~\cite{BedfordFisher1} consider the following averaging integral: 
$$D_d^2(\mu,x) := \lim_{\varepsilon\downarrow 0} \frac{1}{|\ln\varepsilon|}\int_\varepsilon^1 \frac{1}{r} \Theta_d^r(\mu,x)  dr,$$
which they call an \emph{order-two density} of $\mu$ at $x$. It is known (cf.~\cite{Falconer2, PatzschkeZahle1, Zahle3}) that  for a class of sets including strictly self-similar $d$-fractals $D_d^2(\Hdr, x)$ is positive, finite and constant $\Hdr$-a.e. We shall denote this $\Hdr$-a.e. constant as $D_d^2(A)$.

In this paper we examine the limiting case as $s\uparrow d$ of the Riesz potential and energy of a measure $\mu$ by considering the following normalized $d$-energy and $d$-potential: $$\Id(\mu) := \lim_{s\uparrow d} (d-s)I_s(\mu) \qquad \Ud{\mu}(x):= \lim_{s\uparrow d} (d-s)U_s^\mu(x),$$ when they exist. In~\cite{Zahle2}, Z\"ahle provides conditions on a measure $\mu$ for which $D_d^2(\mu,\cdot)$ and $\Ud{\mu}$ agree. (cf.~\cite{Hinz1} for generalizations to other averaging schemes.) We use this result to prove that the limit $\Id(\mu)$ exists for all $\mu\in\MA$, that this normalized energy gives rise to a minimization problem with a unique solution and use this minimization problem to study the behavior of the equilibrium measures $\mu^{s,A}$ as $s\uparrow d$. 

The study of Riesz potentials on fractals is also examined in~\cite{Zahle1, Zahle5} by Z\"ahle in the context of harmonic analysis on fractals. In~\cite{Putinar1}, Putinar considers a different normalization for the Riesz $d$-potential in his work on inverse moment problems.

\subsection{Results}

\begin{theorem}
\label{th:Existence}
Let $A$ be a strictly self-similar $d$-fractal and let $\lambda^d := \Hdr/\Hd(A)$, then

\begin{itemize}

\item[(1)] The limit $\Id(\mu)$ exists for all $\mu\in\MA$ and
$$\Id(\mu)=\left\{ 
\begin{array}{cc}
d D_d^2(A) \int \left(\frac{d\mu}{d\Hdr}\right)^2 d\Hdr & \qquad\text{if }\mu\ll\Hdr,\\ 
\infty & \qquad\text{otherwise.}\end{array}
\right.$$ 

\item[(2)] If $\Id(\mu)<\infty$, then the limit $\Ud{\mu}$ equals $\frac{d\mu}{d\Hdr}$ $\mu$-a.e. and
$$
\Id(\mu) = \int \Ud{\mu} d\mu. 
$$

\item [(3)] $\Id(\lambda^d)<\Id(\nu)$ for all
$\nu\in\MAp\backslash\left\{ \lambda^d\right\} $.
\end{itemize}
\end{theorem}

\begin{theorem}
\label{th:Growth}
Let $A$ be a strictly self-similar $d$-fractal, then there is a constant $K$ depending only on $A$, so that for any $s\in(0,d)$, $\mu^{s,A}(B(x,r)) \le K r^s$ for $\mu^{s,A}$-a.a. $x\in A$ and $r>0$.
\end{theorem}

A bound similar to that in Theorem~\ref{th:Growth} is presented in \cite[Ch. 8]{Mattila1}. This result differs in that the constant $K$ does not depend on $s$. 

\begin{theorem}
\label{th:Convergence}
Let $A$ be a strictly self-similar $d$-fractal and let $\lambda_d := \Hdr/\Hd(A)$, then $\mu^{s,A}$ converges in the weak-star topology on $\MA$ to $\lambda^d$ as $s\uparrow d$.
\end{theorem}

\section{The Existence of a Unique Minimizer of $\Id$}\label{sec:Existence}

A set $A$ is said to be {\em Ahlfors $d$-regular} if there are constants $0<C_1$, $C_2<\infty$ depending only on $A$ so that for all $x$ in $A$ and all $r\in(0,\diam A)$ 
$$C_1r^d < \Hdr(B(x,r)) < C_2r^d.
$$

The proof of Proposition~\ref{prop:ADR} is given by Hutchinson in~\cite[\S 5.3]{Hutchinson1}.

\begin{proposition}
\label{prop:ADR} If $A$ is a strictly self-similar $d$-fractal, then $A$ is Ahlfors $d$-regular.\end{proposition}

The potential $U_s^\mu(x)$ of a finite Borel measure $\mu$ at a point $x$ has the following useful expression in terms of densities: (cf.~\cite{Mattila1}) 
\begin{eqnarray*}
U_s^\mu(x) &:=& \int \Rk{x}{y}{s} d\mu(y) \\
&=& \int_0^\infty \mu(\{ y : |x-y|^{-s} \ge t\}) dt \\
&=& \int_0^\infty \mu(\{ y : |x-y| \le t^{-1/s}\}) dt \\
&=& s\int_0^\infty \frac{\mu(B(x,r))}{r^{s+1}} dr = s\int_0^\infty \Theta_d^r(\mu, x)  \frac{1}{r^{1-(d-s)}} dr,
\end{eqnarray*}
where the second to last equality results from a change of variables replacing $t^{-1/s}$ with $r$. Note that for all $R>0$ $$\lim_{s\uparrow d} (d-s) s\int_R^\infty \Theta_d^r(\mu, x)  \frac{1}{r^{1-(d-s)}} dr = 0.$$ From this we conclude that if $\Ud{\mu}(x)$ exists, then $$\Ud{\mu}(x) = \lim_{s\uparrow d} (d-s) s\int_0^R \Theta_d^r(\mu, x)  \frac{1}{r^{1-(d-s)}} dr,$$ for any $R>0$.

The relationship between the order-two density and the limiting potential is examined by Z\"ahle in the context of stochastic differential equations in~\cite{Zahle2} and also by Hinz, in~\cite{Hinz1}. We include a proof of this relationship from~\cite{Hinz1}. 

\begin{proposition} \label{prop:Order2DensEqNormPot}
Let $\mu$ be a finite Borel measure with support in $\Rp$, $x\in\supp \mu$, $d\in(0,p]$. If $D_d^2(\mu,x)$ exists and is finite, then  $$\Ud{\mu}(x) = dD_d^2(\mu,x).$$
\end{proposition}

\begin{proof}
One may verify that the function $k_\varepsilon(t) := \varepsilon^2 \chi_{(0,1]}(t)  t^{\varepsilon-1} |\log t |$ is an approximate identity in the following sense: If $f:\Rl\to\Rl$ is right continuous at $0$ and is bounded on $(0,1)$, then $$\lim_{\varepsilon\downarrow 0} \int_0^\infty k_\varepsilon(t) f(t) dt = f(0).$$ Define the following function:
$$
f(t) := \left\{
\begin{array}{cc}
\frac{1}{|\log t |}\int_t^1 \frac{1}{r} \Theta_d^r(\mu, x) dr & \text{ when $t>0$} \\
D_d^2(\mu, x) & \text{ when $t=0$} 
\end{array}
\right.
$$
If $D_d^2(\mu, x)$ exists and is finite, then $f$ is right-continuous at $0$ and bounded on $(0,1)$ thus 
\begin{eqnarray*}
D_d^2(\mu, x) &=& \lim_{\varepsilon \downarrow 0} \int_0^\infty k_\varepsilon(t) f(t) dt \\
&=& \lim_{\varepsilon \downarrow 0} \varepsilon^2 \int_0^1 t^{1-\varepsilon} \int_0^1 \frac{\chi_{[t,1]}(r)}{r} \Theta_d^r(\mu, x) dr dt \\
&=& \lim_{\varepsilon \downarrow 0}\varepsilon^2  \int_0^1 \frac{1}{r} \Theta_d^r(\mu, x)\int_1^r t^{\varepsilon -1} dt dr \\
&=& \lim_{\varepsilon \downarrow 0}\varepsilon\int_0^1 \frac{1}{r} \Theta_d^r(\mu, x) r^\varepsilon dr \\
&=&  \lim_{s\uparrow d} (d-s) \int_0^1 \Theta_d^r(\mu, x)  \frac{1}{r^{1-(d-s)}} dr = \frac{1}{d}\Ud{\mu}(x).
\end{eqnarray*}
\end{proof}

\begin{lemma} \label{lemma:measures}
Let $A$ be a strictly self-similar $d$-fractal and let $\mu\in\MA$. If $\mu = \mu^\ll + \mu^\perp$ is the Lebesgue decomposition of $\mu$ with respect to $\Hdr$, then 
\begin{itemize}
\item[1.] $\Ud{\mu^\perp}(x) = \infty$ for $\mu^\perp$-a.a. $x$.
\item[2.] $\Ud{\mu^\ll}(x) = d D_d^2(A) \frac{d\mu}{d\Hdr}(x)$ for $\mu^\ll$-a.a. $x$.
\end{itemize}
\end{lemma}
\begin{proof}
The Radon-Nikod\'ym theorem ensures that for $\mu^\perp$-a.a. $x$, $$\lim_{r\downarrow 0} \frac{\mu^\perp(B(x,r))}{\Hdr(B(x,r))} =\infty.$$ For such an $x$, let $M\in\Rl$ be arbitrary and $R>0$ such that for all $r\in(0,R)$ we have $\mu^\perp(B(x,r))/\Hdr(B(x,r)) > M.$ It then follows that
\begin{eqnarray*}
\liminf_{s\uparrow d} (d-s)s\int_0^\infty \frac{\mu^\perp(B(x,r))}{r^{s+1}} dr &\ge & \left(\inf_{r\in(0,R)} \frac{\mu^\perp(B(x,r))}{\Hdr(B(x,r))}\right) \liminf_{s\uparrow d} (d-s)s\int_0^R \frac{\Hdr(B(x,r))}{r^{s+1}}dr \\
&\ge & M \lim_{s\uparrow d} (d-s)s C_1 \frac{1}{d-s} R^{d-s} = C_1Md,
\end{eqnarray*} 
where $C_1$ is the lower bound from the Ahlfors $d$-regularity of $A$. $M$ is arbitrary, and this proves the first claim.

To prove the second claim we begin with the following equality for an arbitrary $R>0$:
\begin{eqnarray}
\nonumber
& &(d-s)s\int_0^R \frac{\mu^{\ll}(B(x,r))}{r^{s+1}} dr \\
\label{eq:2}
&=& \frac{d\mu^{\ll}}{d\Hdr}(x)(d-s) s\int_0^R \frac{\Hdr(B(x,r))}{r^{s+1}} dr + 
(d-s)s \int_0^R \left(\frac{\mu^{\ll}(B(x,r))}{\Hdr(B(x,r))} - \frac{d\mu^{\ll}}{d\Hdr}(x)\right) \frac{\Hdr(B(x,r))}{r^{s+1}} dr.
\end{eqnarray}
By Proposition~\ref{prop:Order2DensEqNormPot} the limit as $s\uparrow d$ of the first summand in~\eqref{eq:2} is $\frac{d\mu^{\ll}}{d\Hdr}(x)dD_d^2(A)$ for $\Hdr$-a.a. $x$. The absolute value of the limit superior of the second summand in~\eqref{eq:2} is bounded for $\Hdr$-a.a. $x$ by $$\sup_{r\in(0,R)}\left|\frac{\mu^{\ll}(B(x,r))}{\Hdr(B(x,r))} - \frac{d\mu^{\ll}}{d\Hdr}(x)\right|dD_d^2(A),$$ which can be made arbitrarily small by choosing $R$ sufficiently small. Thus the limit as $s\uparrow d$ of~\eqref{eq:2} exists $\Hdr$-a.e. and hence $\Ud{\mu^{\ll}}$ does as well. 
\end{proof}

For a measure $\mu\in\MA$, let   
$$\underline{\Id}(\mu):=\liminf_{s\uparrow d} (d-s) \iint \Rk{x}{y}{s} d\mu(y) d\mu(x).$$

\begin{proposition} 
\label{prop:FiniteEnImpliesL2}
Let $A$ be a strictly self-similar $d$-fractal. If $\underline{\Id}(\mu) <  \infty$ for $\mu\in\MA$, then $\mu \ll \Hdr$ and $\frac{d\mu}{d\Hdr} \in L^2(\Hdr)$.
\end{proposition}

\begin{proof}
Let $\mu\in\MA$ so that $\underline{\Id}(\mu)<\infty$, then by Fatou's lemma $$\int \liminf_{s\uparrow d} (d-s)U_s^\mu d\mu \leq \underline{\Id}(\mu) < \infty.$$ This implies that $\liminf_{s\uparrow d} (d-s)U_s^\mu$ is finite $\mu$-a.e. and, by the first claim in Lemma \ref{lemma:measures}, $\mu \ll \Hdr$. By the second claim in Lemma \ref{lemma:measures} and the previous equation
$$ \int \left(\frac{d\mu}{d\Hdr}\right)^2 d\Hdr =  \int  \left(\frac{d\mu}{d\Hdr}\right) d\mu =  \int \frac{1}{dD_2^d(A)}\Ud{\mu} d\mu  < \infty.$$
\end{proof}

\subsection{Proof of Theorem 1.1}

With the preceding results we may now prove Theorem~\ref{th:Existence}.

\begin{proof}
Let $\mu\in \MA$ so that $\underline{\Id}(\mu) < \infty$, then $\mu\ll\Hdr$ and $d\mu/d\Hdr \in L^2(\Hdr)$. The maximal function of $\mu$ with respect to $\Hdr$ is $$M_{\Hdr}\mu (x) := \sup_{r>0} \frac{\mu(B(x,r))}{\Hdr(B(x,r))} =  \sup_{r>0}\frac{1}{\Hdr(B(x,r))}\int_{B(x,r)} \frac{d\mu}{d\Hdr} d\Hdr.$$ The maximal function is bounded on $L^2(\Hdr)$ and so $M_{\Hdr}\mu\in L^2(\Hdr)$. We shall use this to provide a $\mu$-integrable bound for $(d-s)U_s^\mu$ that is independent of $s$ and appeal to dominated convergence. We begin with the point-wise bound
\begin{eqnarray}
\nonumber
(d-s) \int\Rk{x}{y}{s} d\mu(y) &=& (d-s)s\int_0^\infty \frac{\mu(B(x,r))}{\Hdr(B(x,r))} \frac{\Hdr(B(x,r))}{r^{s+1}} dr \\
\nonumber
&\le& M_{\Hdr}\mu(x) (d-s)s\left[ \int_0^{\diam A} \frac{\Hdr(B(x,r))}{r^{s+1}} dr + 
\int_{\diam A}^\infty  \frac{\Hdr(B(x,r))}{r^{s+1}} dr \right] \\
\label{eq:3}
&\le& M_{\Hdr}\mu(x) \left[ (d-s)s\int_0^{\diam A} \frac{C_2r^d}{r^{s+1}} dr + 
(d-s)s\int_{\diam A}^\infty  \frac{1}{r^{s+1}} dr \right],
\end{eqnarray}
where $C_2$ is the constant in the upper bound of the Ahlfors $d$-regularity of $A$. The quantity in brackets in \eqref{eq:3} may be maximized over $s\in(0,d)$ and we denote this maximum by $K$. Then, by the Cauchy-Schwarz inequality, 
$$\int K M_{\Hdr}\mu\, d\mu <  K\int  \left(M_{\Hdr}\mu\right)\left( \frac{d\mu}{d\Hdr}\right) d\Hdr < K\left\|M_{\Hdr}\mu \right\|_{2,\Hdr}\left\|\frac{d\mu}{d\Hdr} \right\|_{2,\Hdr}<\infty.$$ By dominated convergence the second claim follows. The first claim follows from the second and from Lemma~\ref{lemma:measures} and Proposition~\ref{prop:FiniteEnImpliesL2}.

The final claim of the theorem follows from a straightforward Hilbert space argument. Let $\nu$ denote the finite measure ${dD_d^2(A)}^{-1}\Hdr$. By Proposition~\ref{prop:FiniteEnImpliesL2} the set of measures with finite normalized $d$-energy is identified with the non-negative cone in $L^2(\nu)$ (denoted $L^2(\nu)_+$) via the map $\mu\leftrightarrow d\mu/d\nu$. Under this map we have $\Id(\mu)=\left\| d\mu/d\nu\right\|_{2,\nu}^2$. A measure $\mu$ of finite $d$-energy is a probability measure if  and only if  $\left\|d\mu/d\nu\right\|_{1,\nu} = 1$. We seek a unique non-negative function $f$ that minimizes $\|\cdot\|_{2,\nu}$ subject to the constraint $\|f\|_{1,\nu}=1$. The non-negative constant function $1/\nu(\Rp)$ satisfies the constraint $\|1/\nu(\Rp)\|_{1,\nu}=1$. Let $f\in L^2(\nu)_+$ such that $\|f\|_{1,\nu}=1$ and $\|f\|_{2,\nu} \le \|1/\nu(\Rp) \|_{2,\nu}$, then
$$
\frac{1}{\nu(\Rp)}=\left\|\frac{f}{\nu(\Rp)}\right\|_{1,\nu} = \left\langle f, \frac{1}{\nu(\Rp)} \right\rangle_\nu 
\leq \|f\|_{2,\nu}\left\|\frac{1}{\nu(\Rp)}\right\|_{2,\nu} 
\leq \left\|\frac{1}{\nu(\Rp)}\right\|^2_{2,\nu} = \frac{1}{\nu(\Rp)}.
$$
Thus
$$
\left\langle f, \frac{1}{\nu(\Rp)} \right\rangle_\nu = \|f\|_{2,\nu}\left\|\frac{1}{\nu(\Rp)}\right\|_{2,\nu}.
$$
From the Cauchy-Schwarz inequality $f=1/\nu(\Rp)$ $\nu$-a.e. By the identification above the measure $\lambda^d := \Hdr/\Hd(A)\in\MAp$, uniquely minimizes $\Id$ over $\MAp$.
\end{proof}

\section{The Weak-Star Convergence and Bound on the Growth of $\mu^s$} \label{sec:Convergence}

The proofs of Theorems~\ref{th:Growth} and~\ref{th:Convergence} rely on the following classical results from Potential Theory (cf.~\cite{Landkof1, Fuglede1}). Let $\Es$ denote the set of signed Radon measures with finite total variation such that $\mu\in\Es$ if and only if $I_s(|\mu|) < \infty$. The set $\Es$ is a vector space and when combined with the following positive-definite bilinear form $$I_s(\mu, \nu) := \iint \Rk{x}{y}{s} d\mu(y)d\nu(x),$$ is a pre-Hilbert space. Further, the minimality of the $s$-energy of $\mu^{s,A}$ implies $U_s^{\mu^{s,A}} = I_s(\mu^{s,A})$ $\mu^{s,A}$-a.e.

We shall also use the Principle of Descent:  Let $\{\mu_n\}_{n=1}^\infty\subset\MA$ be a sequence of measures converging in the weak-star topology on $\MA$ to $\psi$ (we shall denote such weak-star convergence with a starred arrow, i.e. $\mu_n\cws \psi$) then for $s\in(0,d)$ $$I_s(\psi) \le \liminf_{n\to\infty} I_s(\mu_n).$$ 

\begin{lemma}
\label{lemma:goesToInfty}
Let $A$ be a compact set for which there is a $C>0$ such that $\Id(\mu)>C$ for all $\mu\in\MAp$, then $$\lim_{s\uparrow d} I_s(\mu^{s,A}) = \infty.$$
\end{lemma}
 
\begin{proof}
Without loss of generality we shall assume that $\diam A \le 1$, then for $0<s<t<d$ and any measure $\mu\in\MA$, $I_s(\mu) \le I_t(\mu)$. For sake of contradiction, assume that the claim does not hold, then there is sequence $\{s_n\}_{n=1}^\infty$ increasing to $d$ so that $$\lim_{n\to\infty} I_{s_n} (\mu^{s_n, A}) = L < \infty.$$ Let $\psi$ be a weak-star cluster point of $\{\mu^{s_n,A}\}_{n=1}^\infty$ (hence a probability measure), and let $\{s_m\}_{m=1}^\infty \subset \{s_n\}_{n=1}^\infty$ so that $\mu^{s_m,A} \cws \psi$. 

For any $s\in(0,d)$ we have 
$$
(d-s)I_s(\psi) \le (d-s)\liminf_{m\to\infty} I_s(\mu^{s_m, A}) \le (d-s)\liminf_{m\to\infty} I_{s_m}(\mu^{s_m, A}) \le (d-s)L.
$$
Letting $s\uparrow d$ implies $\Id(\psi) = 0$, which is a contradiction.
\end{proof}

\begin{lemma}
\label{lemma:AsympDensity}
Let $A$ be a compact set for which there is a $C>0$ such that $\Id(\mu)>C$ for all $\mu\in\MAp$, then $$\lim_{s\uparrow d} \sup_{y\in A} \dist(y, \supp \mu^{s,A}) = 0.$$
\end{lemma}

\begin{proof}
Let $s\in(0,d)$ and $\delta=\sup_{y\in A}\dist(y, \supp \mu^{s,A})$. We consider the possibility that $\delta>0$. Pick $y'\in A$ so that $\dist(y', \supp \mu^{s,A}) > \delta/2$. Let $\nu = \Hd_{A\cap B(y',\delta/4)}/\Hdr(B(y',\delta/4))$. For $\beta\in[0,1]$ we have $(1-\beta)\mu^{s,A} + \beta \nu\in\MAp$. Arguments similar to those used in the proof of Lemma~\ref{lemma:measures} show that $I_s(\nu)<\infty$ for all $s\in(0,d)$.
Define the function 
$$f(\beta) := I_s\left((1-\beta)\mu^{s,A} + \beta \nu\right) = (1-\beta)^2I_s(\mu^{s,A}) + \beta^2 I_s(\nu) + 2\beta(1-\beta)I_s(\mu^{s,A}, \nu). $$ 
Differentiating gives
$$
\frac{1}{2}\frac{df}{d\beta} = \beta\left[I_s(\mu^{s,A} - \nu)\right] - \left[I_s(\mu^{s,A}) - I_s(\mu^{s,A},\nu)\right] 
\qquad \text{and} \qquad \frac{1}{2}\frac{d^2f}{d\beta^2} = \left[I_s(\mu^{s,A} - \nu)\right]. 
$$
Because $I_s(\cdot,\cdot)$ is positive definite, $I_s(\mu^{s,A} - \nu) > 0$. Because $\mu^{s,A}$ is the unique minimizer of $I_s$, $f$ cannot have a minimum for any $\beta >0$, hence $I_s(\mu^{s,A}) - I_s(\mu^{s,A},\nu) \le 0$. We obtain
$$
I_s(\mu^{s,A}) \le I_s(\mu^{s,A},\nu) \le \frac{1}{(\delta/4)^s}, 
\qquad
\text{and hence}
\qquad
\delta \le \frac{4}{I_s(\mu^{s,A})^{1/s}}.
$$
By Lemma \ref{lemma:goesToInfty} $\delta \downarrow 0$ as $s\uparrow d$.
\end{proof}
\subsection{Proof of Theorem 1.2}

The next lemma is straightforward and its proof, which is included for completeness, employs common techniques and ideas presented by e.g. Hutchinson in~\cite{Hutchinson1}. For the rest of the paper we shall order our maps $\{\varphi_1,\ldots,\varphi_N\}$ so that the scaling factors satisfy $L_1\le L_2\le \ldots\le L_N$. 

\begin{lemma}
\label{lemma:Copies}
Let $A$ be a strictly self-similar $d$-fractal then, for each $x\in A$ and $r>0$ there is a subset $A'\subset A$ so that 
\begin{itemize}
\item[1.] $B(x,r) \cap A \subset A'$.
\item[2.] $A' = \varphi(A)$ for some similitude $\varphi$.
\item[3.] $\diam A' < W r$ where $W$ depends only on the set $A$.
\end{itemize}
\end{lemma}

\begin{proof}
Choose $x\in A$ and $r>0$. Let $\tilde{K} = \min_{i\in 1,\ldots, N} \{\dist(\varphi_i(A), A\backslash \varphi_i(A))\}$.  If $r\ge L_1\tilde{K}$, let $A'=A$ and then trivially $A\cap B(x,r) \subset A'$ and $\diam A' < r (2\diam A)/(L_1\tilde{K})$. 

We now consider the case when $r<L_1\tilde{K}$.   Because the images of $A$ under each $\varphi_i$ are disjoint, we may assign to every $y\in A$ a unique infinite sequence $\{j_1, j_2\ldots\}\in\{1,\ldots, N\}^{\mathbb{N}}$ so that $\{y\} =\bigcap_{n=1}^\infty\varphi_{j_n}(\varphi_{j_{n-1}}(\ldots\varphi_{j_1}(A)\ldots))$. If $\{i_1,i_2,\ldots\}$ is the sequence identifying $x$, let $M$ be the smallest natural number so that $L_{i_1}L_{i_2}\ldots L_{i_M}\tilde{K} < r$ (note that $M\ge 2$), then $$r\le L_{i_1}L_{i_2}\ldots L_{i_{M-1}}\tilde{K} < \frac{r}{L_{i_M}}< \frac{r}{L_1}.$$ Let $A' = \varphi_{i_{M-1}}(\varphi_{i_{M-2}}(\ldots\varphi_{i_1}(A)\ldots))$, hence $\diam A' = L_{i_1}L_{i_2}\ldots L_{i_{M-1}}\diam A < r\diam A/(L_1\tilde{K})$. To complete the proof we shall show $B(x,r)\cap A \subset A'$.

Choose $y\in B(x,r)\cap A$. If $y=x$, then $y\in A'$, otherwise let $\{j_1, j_2\ldots\}$ be the sequence identifying $y\in A$ and $m$ the smallest natural number so that $j_m \ne i_m$. We have that $$L_{i_1}L_{i_2}\ldots L_{i_{m-1}}\tilde{K} \le \dist(x,y)\le r \le L_{i_1}L_{i_2}\ldots L_{i_{M-1}}\tilde{K},$$ from which we conclude $m\ge M$ forcing $y\in \varphi_{i_{M-1}}(\varphi_{i_{M-2}}(\ldots\varphi_{i_1}(A)\ldots))=A'$.

The claimed constant $W$ is $(2\diam A)/(L_1\tilde{K})$.
\end{proof}

The remaining proofs will make use of the following fact regarding the behavior of equilibrium measures on scaled sets: If $B'=\varphi(B)$ where $\varphi$ is a similitude with a scale factor of $L$, then for any Borel set $E\subset B'$, $\mu^{s,B'}(E) = \mu^{s,B}(\varphi^{-1}(E))$ and $I_s(\mu^{s,B'}) = L^{-s}I_s(\mu^{s,B})$. This follows from scaling properties of the Riesz kernel.

\begin{proof}[Proof of Theorem~\ref{th:Growth}]
Without loss of generality assume $\diam A \le 1$. Let $x\in A$ and $r\in(0,\diam A/4)$, then 
\begin{equation}
\label{eq:4}
I_s(\mu^{s,A}) = I_s\left(\mu^{s,A}_{B(x,r)} + \mu^{s,A}_{A\backslash B(x,r)}\right) \ge 
I_s\left(\mu^{s,A}_{B(x,r)}\right) + I_s\left(\mu^{s,A}_{A\backslash B(x,r)}\right).
\end{equation} 
By Lemma~\ref{lemma:AsympDensity} there is an $s_0 \in(0,d)$ so that $\mu^{s,A}(A\backslash B(x,\diam A/4))>0$ for all $s\in(s_0,d)$. Note that the choice of $s_0$ depends only on $A$ and not on $x$. First, consider the case $s\in(s_0,d)$. If $\mu^{s,A}(B(x,r)) = 0$, then the claim is trivially proven. Assume $\mu^{s,A}(B(x,r))>0$. We normalize the measures on the right hand side of \eqref{eq:4} to be probability measures and obtain
\begin{equation}
\label{eq:5}
I_s\left(\mu^{s,A}_{B(x,r)}\right) + I_s\left(\mu^{s,A}_{A\backslash B(x,r)}\right) = 
\mu^{s,A}(B(x,r))^2I_s\left(\frac{\mu^{s,A}_{B(x,r)}}{\mu^{s,A}(B(x,r))}\right) + (1-\mu^{s,A}(B(x,r)))^2I_s\left(\frac{\mu^{s,A}_{A\backslash B(x,r)}}{1-\mu^{s,A}(B(x,r))}\right).
\end{equation}
By Lemma~\ref{lemma:Copies} we may find a set $A'\subset A$ so that $B(x,r) \cap A \subset A'$, $\diam A' <Wr$ and $A'$ is a scaling of $A$. The right hand side of~\eqref{eq:5} is bounded below by
\begin{eqnarray}
\nonumber
& &\mu^{s,A}(B(x,r))^2I_s(\mu^{s,A'}) + (1-\mu^{s,A}(B(x,r)))^2I_s(\mu^{s,A})\\
\nonumber
&=&I_s(\mu^{s,A})\left[\mu^{s,A}(B(x,r))^2 \left(\frac{\diam A'}{\diam A}\right)^{-s} + (1-\mu^{s,A}(B(x,r)))^2\right]\\
\label{eq:6}
&>&I_s(\mu^{s,A})\left[\mu^{s,A}(B(x,r))^2 \left(\frac{Wr}{\diam A}\right)^{-s} + (1-\mu^{s,A}(B(x,r)))^2\right]
\end{eqnarray}
Combining~\eqref{eq:4} and~\eqref{eq:6} and dividing by $I_s(\mu^{s,A})$ gives the following:
$$
1\ge\mu^{s,A}(B(x,r))^2 \left(\frac{Wr}{\diam A}\right)^{-s} + 1-2\mu^{s,A}(B(x,r)) + \mu^{s,A}(B(x,r))^2,
$$ hence
$$
2\mu^{s,A}(B(x,r))\ge \mu^{s,A}(B(x,r))^2\left[\left(\frac{Wr}{\diam A}\right)^{-s} + 1\right],
\qquad \text{and thus}
\qquad
\mu^{s,A}(B(x,r)) \le 2\left(\frac{W}{\diam A}\right)^{s} r^s.
$$
Let $K_1$ be the maximum of  $2(W/\diam A)^s$ over $s\in[0,d]$, $K_2$ the maximum of  $(4/\diam A)^s$ over $s\in[0,d]$ and $K_a:=\max\{K_1, K_2\}$, then $\mu^{s,A}(B(x,r)) < K_a r^s$ for all $x\in A$, $r>0$ and $s\in(s_0,d)$. 

For $s\in(0,s_0]$ we have the bound (cf.~\cite[Ch. 8]{Mattila1}) $\mu^{s,A}(B(x,r)) \le U_s^{\mu^{s,A}}(x)r^s = I_s(\mu^{s,A})r^s$ for $\mu^{s,A}$-a.a. $x$. Because $\diam A \le 1$, $I_s(\mu^{s,A}) \le I_{s_0}(\mu^{s_0,A})$ for all $s\in(0, s_0]$. Let $K = \max\{K_a, 2I_{s_0}(\mu^{s_0,A})\}$, then $\mu^{s,A}(B(x,r)) < K r^s$ for $\mu^{s,A}$-a.a. $x\in A$ and $r>0$.
\end{proof}

\subsection{Proof of Theorem 1.3}
\begin{proof}[Proof of Theorem~\ref{th:Convergence}]
Let $f:A\to\Rl$ be continuous. Since $A$ is compact $f$ is uniformly continuous on $A$. Fix $\varepsilon >0$ and let $\delta >0$ so that $f(A\cap B(x,\delta)) \subset (f(x)- \varepsilon,f(x)+ \varepsilon)$ for all $x\in A$. Let $M$ be a natural number high enough so that $L_N^M\diam A < \delta$. 

Let $\alpha$ be a multi-index of length $M$ taking values in $\{1,\ldots, N\}^M$. If $\alpha = (i_1,\ldots,i_M)$, then we denote $\varphi_{i_M}(\varphi_{i_{M-1}}(\ldots(\varphi_{i_1})\ldots))$ by $\phi_\alpha$. Let $\tilde{x}$ be any point in $A$. For any $\nu\in\MAp$ we may write $$\int f d\nu = \sum_{\alpha} \int f d\nu_{\phi_\alpha(A)} = \sum_{\alpha} f(\phi_{\alpha}(\tilde{x}))\nu(\phi_{\alpha(A)}) + \sum_{\alpha} \int\left(f - f(\phi_\alpha(\tilde{x}))\right)d\nu_{\phi_\alpha(A)}.$$ It follows that \begin{equation}\label{eq:7}\left| \int f d\nu - \sum_{\alpha} f(\phi_{\alpha}(\tilde{x}))\nu(\phi_{\alpha(A)})\right| <\varepsilon.\end{equation}

As in the proof of Lemma~\ref{lemma:Copies} let $\tilde{K} = \min_{i\in1,\ldots,N} \{\dist(\varphi_i(A), A\backslash \varphi_i(A)) \}$. If $\alpha$ and $\alpha'$ are different multi-indecies of length $M$, then $\dist(\phi_\alpha(A), \phi_{\alpha'}(A)) \ge L_N^{M-1}\tilde{K}$. By Lemma~\ref{lemma:AsympDensity} there is an $s_0<d$ so that for all $s\in(s_0,d)$ we have $\sup_{y\in A}\dist(y, \supp \mu^{s,A}) < L_N^{M-1}\tilde{K}$. From this we conclude $\mu^{s,A}(\phi_\alpha(A)) > 0$ for any multi-index $\alpha$ of length $M$ and any $s\in(s_0,d)$. For such a choice of $s$ we have $$I_s(\mu^{s,A}) > \sum_{\alpha} I_s\left(\mu^{s,A}_{\phi_\alpha(A)}\right) = \sum_{\alpha} \mu^{s,A}({\phi_\alpha(A)})^2I_s\left(\frac{\mu^{s,A}_{\phi_\alpha(A)}}{\mu^{s,A}({\phi_\alpha(A)})}\right) \ge \sum_{\alpha} \mu^{s,A}({\phi_\alpha(A)})^2I_s\left(\mu^{s,\phi_\alpha(A)}\right).$$ We shall use the notation $L_\alpha$ to denote $L_{i_1}L_{i_2} \ldots L_{i_M}$. By appealing to the scaling properties of the Riesz energy, the above becomes 
$$
I_s(\mu^{s,A}) > \sum_{\alpha} \mu^{s,A}({\phi_\alpha(A)})^2L_\alpha^{d-s}\frac{ I_s(\mu^{s,A})}{L_\alpha^d}.
$$  
Let $\psi$ be any weak-star cluster point of $\mu^{s,A}$ as $s\uparrow d$ and let $\{s_n\}_{n=1}^\infty\uparrow d$ be a sequence so that $\mu^{s_n,A}\cws \psi$ and hence so that $(\mu^{s_n,A}(\phi_\alpha(A)))_\alpha$ converges in $[0,1]^{N^M}$, then 
$$
1 = \lim_{n\to\infty}\frac{1}{(L_1^M)^{d-s_n}} \ge \lim_{n\to\infty}\sum_{\alpha} \frac{\mu^{s_n,A}(\phi_\alpha(A))^2}{L_\alpha^d} = \sum_{\alpha} \frac{[\lim_{n\to\infty}\mu^{s_n,A}(\phi_\alpha(A))]^2}{L_\alpha^d}.
$$ 
We then have that 
$$1 = \sum_{\alpha} \lim_{n\to\infty}\mu^{s_n,A}(\phi_\alpha(A)) = \sum_{\alpha} \frac{\lim_{n\to\infty}\mu^{s_n,A}(\phi_\alpha(A)}{\sqrt{L_\alpha^d}} \sqrt{L_\alpha^d} \le \sqrt{\sum_{\alpha} \frac{[\lim_{n\to\infty}\mu^{s_n,A}(\phi_\alpha(A))]^2}{L_\alpha^d}}\sqrt{\sum_{\alpha} L_\alpha^d}= 1.
$$ 
Note that the sum over $\alpha$ of $L_\alpha^d$ is one because the sum over $i\in1,\ldots N$ of $L_i^d$ is one. From this we conclude 
$$
\lim_{n\to\infty} \mu^{s_n,A}(\phi_\alpha(A))=L_\alpha^d
$$
for every multi-index $\alpha$ of length $M$. Because $\lambda^d(\phi_\alpha(A)) = L_\alpha^d$, we have that $$\lim_{n\to\infty} \sum_{\alpha} f(\phi_{\alpha}(\tilde{x}))\mu^{s_n,A}(\phi_\alpha(A)) = \sum_{\alpha} f(\phi_{\alpha}(\tilde{x}))\lambda^d(\phi_\alpha(A)),$$ and so $$\left| \lim_{n\to\infty} \int f d\mu^{s_n,A} - \int f d\lambda^d\right| < 2\varepsilon.$$ The choice of $\varepsilon$ in \eqref{eq:7} was arbitrary as was the choice of the continuous function $f$ and so $\lambda^d=\psi$ for any weak-star cluster point $\psi$, and hence $\mu^{s,A}\cws\lambda^d$ as $s\uparrow d$. 
\end{proof}

\bibliography{References}
\bibliographystyle{abbrv} 

\end{document}